\documentclass[11pt]{amsart}

\usepackage{amsmath}
\usepackage{amsfonts}
\usepackage{amssymb}
\usepackage{amscd}
\usepackage{amsthm}
\usepackage{framed}
\usepackage{fullpage}
\usepackage{graphicx}
\usepackage{latexsym}
\usepackage[numbers]{natbib}

\usepackage{hyperref}
\hypersetup{colorlinks,linkcolor={red},citecolor={blue},urlcolor={blue}}

\newtheorem{theorem}{Theorem}[section]
\newtheorem{lemma}[theorem]{Lemma}
\newtheorem{proposition}[theorem]{Proposition}
\newtheorem{corollary}[theorem]{Corollary}

\theoremstyle{definition}
\newtheorem{definition}[theorem]{Definition}

\numberwithin{equation}{section}

\newcommand{\CC}{\mathbb C}
\newcommand{\HH}{\mathbb H}

\newcommand{\NN}{\mathbb N}

\newcommand{\cH}{\mathcal H}

\newcommand{\QQ}{\mathbb Q}

\newcommand{\ZZ}{\mathbb Z}

\newcommand{\SL}{\mathop{\mathrm {SL}}\nolimits}

\newcommand{\Sp}{\mathop{\mathrm {Sp}}\nolimits}
\newcommand{\Orth}{\mathop{\null\mathrm {O}}\nolimits}

\newcommand{\latt}[1]{{\langle{#1}\rangle}}

\def\dim{\operatorname{dim}}

\def\w{\operatorname{w}}

\newenvironment{psmallmatrix}
  {\left(\begin{smallmatrix}}
{\end{smallmatrix}\right)}

\begin{document}

\title[Sums of four polygonal numbers: precise formulas]{Sums of four polygonal numbers: precise formulas} 

\author{Jialin Li}

\address{School of Mathematics and Statistics, Wuhan University, Wuhan 430072, Hubei, China}

\email{jlli.math@whu.edu.cn}

\author{Haowu Wang}

\address{School of Mathematics and Statistics, Wuhan University, Wuhan 430072, Hubei, China}

\email{haowu.wangmath@whu.edu.cn}

\subjclass[2020]{11E25, 11F50, 11F30, 11D85}

\date{\today}

\keywords{Polygonal numbers, Lagrange's four-square theorem, Jacobi forms, Fourier coefficients of automorphic forms, Representations of integers}

\begin{abstract} 
In this paper we give unified formulas for the numbers of representations of positive integers as sums of four generalized $m$-gonal numbers, and as restricted sums of four squares under a linear condition, respectively. These formulas are given as $\mathbb{Z}$-linear combinations of Hurwitz class numbers. As applications, we prove several Zhi-Wei Sun's conjectures. As by-products, we obtain formulas for expressing the Fourier coefficients of $\vartheta(\tau,z)^4$, $\eta(\tau)^{12}$, $\eta(\tau)^4$ and $\eta(\tau)^8\eta(2\tau)^8$ in terms of Hurwitz class numbers, respectively. The proof is based on the theory of Jacobi forms.
\end{abstract}

\maketitle

\section{Introduction}
Let $m\geq 3$ be an integer and $\ell \in \NN=\{0,1,2,...\}$. The $\ell$-th $m$-gonal number
\begin{equation}
p_m(\ell)=\frac{(m-2)\ell^2-(m-4)\ell}{2}     
\end{equation}
counts the number of dots arranged in the shape of a regular polygon with $m$ sides having $\ell$ dots on each side. In particular, $p_3(\ell)$ and $p_4(\ell)$ give the triangular numbers and squares, respectively. The study of representing positive integers as sums of polygonal numbers has a long and storied history (see e.g. \cite[pp. 3-35]{Nat96}). Fermat conjectured in 1638 that every positive integer can be written as the sum of at most $m$ polygonal numbers of order $m$. This claim was proved by Lagrange in the case $m=4$, Gauss in the case $m=3$, and Cauchy in the case $m\geq 5$. Later, Legendre proved that every integer $N\geq 28m^3$ is the sum of four polygonal numbers of odd order $m+2$, and every odd integer $N\geq 7m^3$ is the sum of four polygonal numbers of even order $m+2$.  

One is interested not only in the existence of representations of integers as sums of polygonal numbers, but also in precise formulas for the numbers of representations. To simplify the problem, we consider generalized polygonal numbers, namely $p_m(\ell)$ with $\ell \in \ZZ$. For any $N\in \NN$ we define the representative number
\begin{equation}
R_{m,t}(N):= \# \big\{ (\ell_1, \ell_2,..., \ell_t)\in \ZZ^t : p_m(\ell_1) + p_m(\ell_2) + \cdots + p_m(\ell_t) = N \big\}.     
\end{equation}
For any positive integer $N$, Jacobi showed that
\begin{equation}\label{eq:Jacobi}
R_{4,4}(N)= 8\sum_{d|N,\; 4\nmid d} d. 
\end{equation}
Formulas like \eqref{eq:Jacobi} are rare because the proofs usually use the theory of elliptic modular forms and there are obstructions from cusp forms. In \cite{GW18}, Gritsenko and the second named author employed the theory of Jacobi forms to give a precise formula for $R_{m,8}(N)$ for any $m\geq 3$ and any $N\in \NN$, which involves Cohen's numbers $H(3,-)$. As far as we know, this is the only exact formula for $R_{m,t}(N)$ for large $m$. 

In this paper we further establish an explicit and unified formula for $R_{m,4}(N)$.  

\begin{theorem}\label{th:4-polygonal}
For any integer $m\geq 1$ and integer $N\geq 0$, 
\begin{equation}\label{eq:R4}
R_{m+2,4}(N)=\sum_{\substack{n\in \NN, \; r\in \ZZ \\ mn/2-r+2-m/2=N \\ 2\nmid n,\; 2|r \\ 4n-r^2\geq 0}} 12H^{(2)}(4n-r^2) \; + \sum_{\substack{n\in \NN, \; r\in \ZZ \\ mn/2-r+2-m/2=N \\ 2\nmid n,\; 2\nmid r \\ 4n-r^2\geq 0}} 12H(4n-r^2),
\end{equation}
where $H(D)$ is the Hurwitz class number and $H^{(2)}(D)$ is the modified Hurwitz class number. 
\end{theorem}

Recall that $H(0)=-1/12$ and $H(D)$ for $D>0$ counts the number of $\SL_2(\ZZ)$-equivalence classes of binary integral positive definite quadratic forms of discriminant $-D$ with multiplicities. The modified Hurwitz class number was introduced by \cite{BSZ19}. If $D\equiv 0, 3 \; (\mathrm{mod}\; 4)$, then  
\begin{align*}
H^{(2)}(D)&=H\big(D/4^a\big)\left( 1 - \left( \frac{-D/4^a}{2} \right) \right)  \\
&=H(4D)-2H(D),
\end{align*}
where $a$ is the largest integer such that $D/4^a \equiv 0, 3 \; (\mathrm{mod}\; 4)$. If $D\equiv 1, 2 \; (\mathrm{mod}\; 4)$, then $H^{(2)}(D)=0$. 

By verifying the non-emptiness of the second sum in \eqref{eq:R4}, we derive the following corollary. This in some sense improves the lower bound of Legendre's theorem. 

\begin{corollary}\label{cor:4-polygonal}
Let $N\geq m^3/8+m/2+2$ be an integer. If $m$ or $N$ is odd, then $R_{m+2,4}(N)>0$. 
\end{corollary}

The case $m=6$ of Theorem \ref{th:4-polygonal} confirms \cite[Conjecture 3.1]{Sun16} for odd integers. 

\begin{corollary}\label{cor:r=1}
Let $r(N)$ be the representative number of $N\in\NN$ as four unordered generalized octagonal numbers. When $N$ is odd, $r(N)=1$ if and only if $N\in\{1,3,5,7,9,13\}$.    
\end{corollary}

When both $m$ and $N$ are even, the second sum of \eqref{eq:R4} is empty. In this case, it is still possible to establish sufficient conditions for $R_{m+2,4}(N)>0$ (comparing with \cite[Theorem 1.1]{MS17}) and prove \cite[Conjecture 3.1]{Sun16} for even integers. 

Theorem \ref{th:4-polygonal} is proved by the theory of Jacobi forms on congruence subgroups. Jacobi forms were first introduced by Eichler and Zagier \cite{EZ85}. They are holomorphic functions of two variables $(\tau,z)\in \HH \times \CC$, which are modular in $\tau$ and quasi periodic in $z$, and satisfy some growth conditions at cusps. It is known (see e.g. \cite[Example 1.5]{GN98}) that the odd Jacobi theta series
\begin{equation}\label{eq:theta}
\begin{split}
\vartheta(\tau,z)=&q^{1/8}\zeta^{1/2}\sum_{n\in \ZZ} (-1)^nq^{(n+1)n/2}\zeta^n \\ =&
-q^{1/8}\zeta^{-1/2}\prod_{n \geq 1} (1-q^{n-1}\zeta)(1-q^{n}\zeta^{-1})(1-q^n)    
\end{split}    
\end{equation}
defines a holomorphic Jacobi form of weight $1/2$ and index $1/2$ with a multiplier system of order $8$ on $\SL_2(\ZZ)$, where $q=e^{2\pi i\tau}$ and $\zeta=e^{2\pi iz}$. We observe that $\vartheta(2\tau,2z)^4$ is a holomorphic Jacobi form of weight $2$, index $4$ and trivial character on $\Gamma_0(4)$. Boylan, Skoruppa and Zhou \cite{BSZ19} proved that 
\begin{equation}
\cH(\tau,z)=\sum_{\substack{n\in \NN, \; r\in\ZZ\\ 4n-r^2\geq 0}} 12H^{(2)}(4n-r^2) q^n\zeta^r 
\end{equation}
generates the $\CC$-vector space of holomorphic Jacobi form of weight $2$ and index $1$ on $\Gamma_0(2)$. By means of the theta decomposition of $\cH(\tau,z)$ and the Hecke operators, we construct a basis of the $3$-dimensional space of holomorphic Jacobi forms of weight $2$ and index $4$ on $\Gamma_0(4)$. This basis has Fourier expansions involving the Hurwitz class numbers $H(-)$. We then deduce the following expression of $\vartheta(2\tau,2z)^4$. 

\begin{theorem}\label{th:power-4}
\begin{equation}\label{eq:power-4}
\vartheta(2\tau,2z)^4 =\sum_{\substack{n\in \NN \\ 2\nmid n}} \left( \sum_{\substack{r\in \ZZ,\; 2|r \\ 4 n-r^2 \geq 0}} 12H^{(2)}(4n-r^2)\zeta^{2r}-\sum_{\substack{r\in\ZZ,\; 2\nmid r \\ 4n-r^2 \geq 0}} 12H(4n-r^2)\zeta^{2r} \right)  q^n.
\end{equation}    
\end{theorem}
Combining \eqref{eq:power-4} and the generating function
\begin{equation}\label{eq:generating}
 \vartheta\big( m\tau, 1/2-\tau \big)^4 = q^{m/2-2} \sum_{N=0}^\infty R_{m+2,4}(N) q^N,
\end{equation}
we complete the proof of Theorem \ref{th:4-polygonal}. 

We further present two applications of Theorem \ref{th:power-4}. A certain specialization of \eqref{eq:power-4} yields

\begin{proposition}\label{prop:eta-4}
\begin{equation}
\prod_{n=1}^{\infty}(1-q^n)^4=\sum_{N=0}^{\infty}\left( \sum_{\substack{n\in \NN, \; r\in \ZZ \\ 2\nmid n, \; 2\mid r \\ 3n+4r+5=2N \\ 4n-r^2\geq 0 }}12H^{(2)}(4n-r^2)-\sum_{\substack{n\in \NN, \; r\in \ZZ \\ 2\nmid n, \; 2\nmid r \\ 3n+4r+5=2N \\ 4n-r^2\geq 0 }}12H(4n-r^2)\right)q^N.  
\end{equation}    
\end{proposition}

Taking the fourth derivative of \eqref{eq:power-4} with respect to $z$ and then setting $z=0$, we find

\begin{proposition}\label{prop:eta-12}
\begin{equation}
q\prod_{n=1}^\infty (1-q^{2n})^{12} =    \frac{1}{2}\sum_{\substack{n\in \NN \\ 2\nmid n}} \left( \sum_{\substack{r\in \ZZ, \; 2|r \\ 4n-r^2\geq 0}}r^4 H^{(2)}(4n-r^2)-\sum_{\substack{ r\in \ZZ, \; 2\nmid r \\ 4n-r^2\geq 0}}r^4 H(4n-r^2) \right)  q^n.   
\end{equation}    
\end{proposition}

Similarly, by analyzing holomorphic Jacobi forms of weight $2$ and index $9$ on $\Gamma_0(2)$, we determine the following expression. 

\begin{proposition}\label{prop:eta-1-2-8}
Let $\tau_2(n)$ denote the coefficient of $q^n$ in the Fourier expansion of the new form $q\prod_{n=1}^\infty (1-q^n)^8 (1-q^{2n})^8$. Then
\begin{equation}
\tau_2(n) =  \frac{1}{23040} \left( 13 \sum_{\substack{r\in \ZZ, \; 3|r \\ 36n\geq r^{2} }} r^6 H^{(2)}\left(\frac{36n-r^2}{9}\right) - \sum_{\substack{r\in \ZZ \\ 36n\geq r^{2}}} \sum_{\substack{a>0 \\ a \mid (n,r,9) }} a r^6 H^{(2)} \left(\frac{36n-r^2}{a^2}\right) \right).  
\end{equation}    
\end{proposition}

It is well-known that the Fourier expansion of $\Delta(\tau)=q\prod_{n=1}^\infty (1-q^n)^{24}$ has an expression involving the Hurwitz class numbers by the Eichler--Selberg trace formula \cite{Sel56}. This approach may also yield similar formulas for the Fourier expansions of the three cusp forms above. It would be interesting to compare the formulas obtained by the two different methods. 

We now extend the Jacobi forms approach to calculate the number of representations of integers as restricted sums of four squares under a linear condition. This is closely related to Sun's work \cite{Sun17, SS18, Sun19} to refine Lagrange's four-square theorem. Let $\Orth(\ZZ^4)$ be the group generated by all permutations and all sign changes acting on $\ZZ^4$. Let $m\in \NN$ and $\mathcal{S}_m$ denote the set of all vectors of (square) length $m$, namely
$$
\mathcal{S}_m = \{ (a,b,c,d) \in \ZZ^4: a^2+b^2+c^2+d^2=m\}. 
$$
We fix an integer $r$. For any $v=(a,b,c,d)\in \mathcal{S}_m$ and any positive integer $n$, we define
\begin{equation}
\rho_v(n,r) = \# \{ (x, y, z, w) \in \ZZ^4 : x^2+y^2+z^2+w^2=n, \; ax+by+cz+dw=r \}.    
\end{equation}
The group $\Orth(\ZZ^4)$ acts on $\mathcal{S}_m$, and the number $\rho_v(n,r)$ depends only on the $\Orth(\ZZ^4)$-orbit of $v$. We label the orbits of $\mathcal{S}_m$ under the action of $\Orth(\ZZ^4)$ by $\Orth(\ZZ^4)\cdot v_i$ for $1\leq i \leq t_m$. Let $\delta_i=| \Orth(\ZZ^4)\cdot v_i|$ denote the length of the orbit. We prove the following weighted formula. 

\begin{theorem}\label{th:four-squares}
For any $m\in \NN$, $r\in \ZZ$ and any positive integer $n$,  
\begin{align*}
\rho(n,r,m):=&\sum_{j=1}^{t_m} \delta_j \cdot \rho_{v_j}(n,r) \\
=& \begin{cases}
96 \left(  \sum_{\substack{ d>0,\; 2\nmid d \\ d\mid(n,r,m)}} dH^{(2)}\Big(\frac{4(nm-r^2)}{d^2}\Big)+2\sum_{\substack{ d>0,\; 2\nmid d \\ d\mid(n,r,m)}}dH^{(2)}\Big(\frac{nm-r^2}{d^2}\Big) \right),& 2\mid (n,m), \\
96\sum_{\substack{ d>0,\; 2\nmid d \\ d\mid(n,r,m)}} dH^{(2)}\Big(\frac{4(nm-r^2)}{d^2}\Big), & 2\nmid (n,m).
\end{cases}
\end{align*}
\end{theorem}

As direct consequences, $\rho(N,0,0)$ recovers Jacobi's formula for $R_{4,4}(N)$ and
\begin{equation}\label{eq:Gauss}
R_{4,3}(N)= \frac{1}{8}\rho(N,0,1) = 12 H^{(2)}(4N),   
\end{equation}
which counts the representative number of $n$ as the sum of three squares. In addition, 
$$
\#\big\{(x,y,z)\in\ZZ^3: N=x^2+y^2+2z^2 \big\} = \frac{1}{24}\rho(N,0,2) = 4H^{(2)}(8N) + 8H^{(2)}(2N),
$$
which follows that the quadratic form $x^2+y^2+2z^2$ represents a positive integer $N$ if and only if $N$ is not of the form $2^{2a+1}(8b+7)$ for any $a,b\in\NN$. This recovers \cite[Theorem 2]{Dur16}.

We derive the following lemma from Gauss--Legendre's three-square theorem.

\begin{lemma}\label{lem:positive}
Let $D$ be a positive integer. Then there is $s\in \NN$ such that $H^{(2)}\big(4(D-16^s)\big)>0$ if and only if $D\not\in \{2^{4k+3} : k\in\NN \}$. As a consequence, given positive integers $m$ and $n$, there is $s\in \NN$ such that $\rho(n,4^s,m)>0$ if and only if $nm\not\in \{2^{4k+3} : k\in\NN \}$. 
\end{lemma} 

The class number $t_m$ is one for some particular integers $m$. In this case, we conclude the following result from Lemma \ref{lem:positive}.   

\begin{corollary}\label{coro:sun}
\noindent
\begin{enumerate}
\item Let $(a,b,c,d)$ be one of the following vectors:
\begin{align*}
&(1,1,1,0), \; (2,1,0,0),\; (2,1,1,0), \; (2,1,1,1), \;  (3,1,1,0), \; (3,2,1,0), \\
&(3,2,1,1), \; (3,3,2,1), \; (4,2,2,0), \; (6,4,2,0), \; (2,2,1,0). 
\end{align*}
For any positive integer $n$, there exist integers $x$, $y$, $z$, $w$ and $s\geq 0$ such that $n=x^2+y^2+z^2+w^2$ and $ax+by+cz+dw=4^s$. 
\item Let $(a,b,c,d)$ be one of the following vectors:
$$
(1,0,0,0), \; (1,1,0,0),\; (2,2,0,0).
$$
Let $m=a^2+b^2+c^2+d^2$ and $n$ be a positive integer. There exist integers $x$, $y$, $z$, $w$ and $s\geq 0$ such that $n=x^2+y^2+z^2+w^2$ and $ax+by+cz+dw=4^s$ if and only if $n\neq \frac{2^{4k+3}}{m}$ for any integer $k\geq 0$. 
\end{enumerate}
\end{corollary}

This corollary recovers or improves some of Sun's previous results. The cases $(1,1,1,0)$, $(2,1,0,0)$, $(2,1,1,0)$ and $(2,2,1,0)$ were proved in \cite{Sun19}. The case $(3,2,1,0)$ was proved in \cite[Theorem 1.8]{Sun20add}. The other cases improve some results in \cite{Sun17, SS18, Sun19}. 

We explain the proof of Theorem \ref{th:four-squares}. Recently, Li and Zhou \cite{LZ24} established a similar formula for the representative number associated with the maximal order in the quaternion algebra $(-1,-1)_{\QQ}$. They proved the formula by identifying the corresponding generating function with the unique Siegel Eisenstein series of weight $2$, level $2$ and degree $2$. Inspired by their proof, we show that the generating function of $\rho(n,r,m)$ defines a Siegel modular form of weight $2$, level $4$ and degree $2$, and construct a basis of the associated $3$-dimensional space of Siegel modular forms by lifting the Jacobi forms of weight $2$ constructed above. We then derive the desired formula for $\rho(n,r,m)$ involving the Hurwitz class numbers. We remark that \cite[Theorem 1.2]{LZ24} can be deduced by taking $m=35$ in Theorem \ref{th:four-squares}.

When $v=(a,b,c,d)$ takes other vectors of $\ZZ^4$, the claim of Corollary \ref{coro:sun} may still hold. Sun proposed many conjectures for such vectors, such as $v=(3,1,0,0)$ and $(5,3,1,0)$. The two cases were proved in \cite{SW21} and \cite{MT21, MRT22}, respectively.  Some similar results were proved in \cite{SZ24} recently.  

The paper is organized as follows. In Section \ref{sec:Jacobi} we briefly introduce Jacobi forms on congruence subgroups. In Section \ref{sec:4-polygonal} we present proofs of Theorems \ref{th:4-polygonal}, \ref{th:power-4}, Corollaries \ref{cor:4-polygonal}--\ref{cor:r=1} and Propositions \ref{prop:eta-4}--\ref{prop:eta-1-2-8}. In Section \ref{sec:4-square} we prove Theorem \ref{th:four-squares} and Corollary \ref{coro:sun}.  
In the appendix, we establish similar formulas for the representative number of positive integers as sums of four generalized $m$-gonal numbers with coefficients $(1,1,3,3)$, and for the Fourier coefficients of $\vartheta(\tau,z)^2\vartheta(3\tau,3z)^2$, $\eta(\tau)^2\eta(3\tau)^2$ and $\eta(\tau)^6\eta(3\tau)^6$, respectively. 

\bigskip
\noindent
\textbf{Acknowledgements} 
H. Wang thanks Prof. Zhi-Wei Sun for his valuable comments and helpful discussions.  The authors thank the referee for careful reading and valuable comments.

\section{Jacobi forms on congruence subgroups}\label{sec:Jacobi}
The theory of Jacobi forms was developed by Eichler and Zagier \cite{EZ85}. In this section, we review some of the basics of Jacobi forms that are necessary to prove the main results of the introduction. 

\begin{definition}
Let $k\in\ZZ$, $t\in \NN$ and $\Gamma$ be a congruence subgroup of $\SL_2(\ZZ)$. A holomorphic function $\phi: \HH \times \CC \to \CC$  is called a \textit{weak Jacobi form} of weight $k$ and index $t$ on $\Gamma$ if it satisfies 
$$
\Big(\phi|_{k,t}A\Big)(\tau,z):=(c\tau+d)^{-k}\exp\left( -\frac{2\pi itcz^2}{c\tau+d} \right) \phi\left( \frac{a\tau+b}{c\tau+d}, \frac{z}{c\tau+d} \right) = \phi(\tau,z)
$$
for any $A=\begin{psmallmatrix}
    a & b \\ c & d
\end{psmallmatrix} \in \Gamma$ and 
$$
\exp\big(2\pi it(x^2\tau + 2xz)\big)\phi(\tau, z+x\tau+y) = \phi(\tau,z)
$$
for any $x,y\in \ZZ$, and if for every $\gamma \in \SL_2(\ZZ)$ the function $\phi|_{k,t}\gamma$ has the Fourier expansion
\begin{equation}
\Big(\phi|_{k,m}\gamma\Big)(\tau,z) = \sum_{n\in \frac{1}{h}\NN, \; r\in \ZZ} f_\gamma(n,r)q^n \zeta^r, \quad q=e^{2\pi i\tau}, \; \zeta=e^{2\pi iz},   
\end{equation}
for some positive integer $h$. If $\phi$ further satisfies that $f_\gamma(n,r)=0$ whenever $4nt<r^2$ for any $\gamma\in \SL_2(\ZZ)$, then it is called a \textit{holomorphic Jacobi form}. We denote the $\CC$-vector space of weak (resp. holomorphic) Jacobi forms of weight $k$ and index $t$ on $\Gamma$ by $J_{k,t}^{\w}(\Gamma)$ (resp. $J_{k,t}(\Gamma)$). 
\end{definition}

We focus on the case $\Gamma=\Gamma_0(N)$ in this paper.  We first recall some standard Hecke operators (see e.g. \cite{EZ85} and \cite[Section 2]{CG11}). 

\begin{lemma}\label{lem:U}
Let $\phi(\tau,z)\in J_{k,t}(\Gamma_0(N))$. Then $\phi(\tau,lz)\in J_{k,tl^2}(\Gamma_0(N))$ and $\phi(l\tau,lz)\in J_{k,tl}(\Gamma_0(Nl))$ for any positive integer $l$. 
\end{lemma}

\begin{lemma}\label{lem:V}
Let $\phi(\tau,z)=\sum_{n,r\in \ZZ}f(n,r)q^n\zeta^r\in J_{k,t}(\Gamma_0(N))$.  For any positive integer $l$, we have
$$
\Big(\phi|_{k,t}V_l\Big)(\tau,z):=\sum_{n,r\in \ZZ} \sum_{\substack{a\geq 1\\ a|(n,r,l) \\ (a,N)=1}} a^{k-1}f(nl/a^2,\, r/a)q^n\zeta^r \in J_{k,tl}(\Gamma_0(N)). 
$$
\end{lemma}

The dimension formula of $J_{2,t}(\Gamma_0(N))$ is unknown in general. We estimate the dimension by the structure of the space of weak Jacobi forms in some specific cases. By \cite{EZ85, Wir92} or \cite[Remark 4.2]{Wan21}, the algebra  
$$
J_{2*,*}^{\w}(\Gamma_0(N)):=\bigoplus_{t\in\NN} \bigoplus_{k\in\ZZ} J_{2k,t}^{\w}(\Gamma_0(N))
$$
is freely generated by $\phi_{-2,1}$ and $\phi_{0,1}$ over $M_{*}(\Gamma_0(N))$, where $\phi_{-2,1}$ (resp. $\phi_{0,1}$) is a basis of $J_{-2,1}^{\w}(\SL_2(\ZZ))$ (resp. $J_{0,1}^{\w}(\SL_2(\ZZ))$) introduced in \cite[Theorem 9.3]{EZ85}), and $M_{*}(\Gamma_0(N))$ is the ring of holomorphic modular forms on $\Gamma_0(N)$. Thus it is easy to find an explicit basis of $J_{2k,t}^{\w}(\Gamma_0(N))$. We characterize $J_{2k,t}(\Gamma_0(N))$ as a specific subspace of $J_{2k,t}^{\w}(\Gamma_0(N))$. Let $\phi\in J_{2k,t}^{\w}(\Gamma_0(N))$.  At each cusp $\gamma \cdot \infty$, the Fourier coefficients $f_\gamma(n,r)$ depend only on the class of $r$ modulo $t$ and on the hyperbolic norm $4nt-r^2$. Consequently, $\phi$ is a holomorphic Jacobi form if and only if at each cusp $\gamma \cdot \infty$ its Fourier coefficients $f_{\gamma}(n,r)$ vanish for all $0\leq r \leq m$ and all $0\leq n < \frac{r^2}{4t}$. Therefore, a linear combination of the basis of $J_{2k,t}^{\w}(\Gamma_0(N))$ lies in the subspace $J_{2k,t}(\Gamma_0(N))$ if and only if the coefficients of this linear combination give a solution to a certain system of linear equations. By this algorithm, we determine the dimensions in Tables \ref{tab:dim-2} and \ref{tab:dim-4}.

\begin{table}[ht]
\caption{The dimension of $J_{2,t}(\Gamma_0(2))$ for $t\leq 25$}
\label{tab:dim-2}
\renewcommand\arraystretch{1.5}
\[
\begin{array}{c|ccccccccccccc}
\hline 
t &  1 & 2 & 3 & 4 &5 & 6&7&8&9&10&11&12&13 \\ 
\hline 
\dim &  1&1&1&1&1&1&1&1&2&1&1&1&2\\
\hline
\hline 
t&14&15&16&17&18&19&20&21&22&23&24&25 \\
\hline
\dim & 1&1&1&2&2&2&1&2&1&1&1&5 \\
\hline
\end{array} 
\]
\bigskip
\end{table}

\begin{table}[ht]
\caption{The dimension of $J_{2,t}(\Gamma_0(4))$ for $t\leq 9$}
\label{tab:dim-4}
\renewcommand\arraystretch{1.5}
\[
\begin{array}{c|ccccccccc}
\hline 
t &  1 & 2 & 3 & 4 &5 & 6&7&8&9 \\ 
\hline 
\dim &  2&2&2&3&2&2&2&2&4 \\
\hline
\end{array} 
\]
\bigskip
\end{table}

\section{A proof of Theorem \ref{th:power-4} and applications}\label{sec:4-polygonal}
In this section, we first prove Theorem \ref{th:power-4}, and then derive some relevant results presented in the introduction from this theorem. 

\subsection{A proof of Theorem \ref{th:power-4}} 
By \cite[Lemma 1.2 and Example 1.5]{GN98}, it is easy to verify that $\vartheta(2\tau,2z)^4$ is a holomorphic Jacobi form of weight $2$ and index $4$ on $\Gamma_0(4)$. It also defines a Jacobi form of quadratic character on $\Gamma_0(2)$. We know from Table \ref{tab:dim-4} that $\dim J_{2,4}(\Gamma_0(4))=3$. We now construct a basis of $J_{2,4}(\Gamma_0(4))$ whose Fourier expansions involve the Hurwitz class numbers. 

As mentioned in the introduction, Boylan, Skoruppa and Zhou \cite{BSZ19} proved that 
$$
\cH(\tau,z)=\sum_{\substack{n\in \NN, \; r\in\ZZ\\ 4n-r^2\geq 0}} 12H^{(2)}(4n-r^2) q^n\zeta^r \in J_{2,1}(\Gamma_0(2)). 
$$
From the theta decomposition
$$
\cH(\tau,z)=h_0(\tau)\sum_{\substack{r\in\ZZ\\ 2|r}} q^{r^2/4}\zeta^r +  h_1(\tau)\sum_{\substack{r\in\ZZ\\ 2\nmid r}} q^{r^2/4}\zeta^r 
=: \cH_0(\tau,z)+\cH_{1}(\tau,z),
$$
we conclude that both $\cH_0$ and $\cH_1$ are holomorphic Jacobi forms of weight $2$ and index $1$ on $\Gamma_0(4)$, since $h_0$, $h_1$ and the two Jacobi theta functions are all modular under $\Gamma_0(4)$ (see \cite[Section 5]{EZ85}). By definition, when $r$ is odd, $H^{(2)}(4n-r^2)$ equals $2H(4n-r^2)$ if $n$ is odd, and equals $0$ otherwise. Thus we have the following Fourier expansions
\begin{align}
\cH_0(\tau, z)=& \sum_{\substack{n,\, r\in \ZZ \\ 2|r \\
					4 n\geq r^2}}12 H^{(2)}\left(4 n-r^2\right) q^n \zeta^r 
			=1+(\zeta^{\pm 2}+6)q+ (6\zeta^{\pm 2}+12)q^2 + O(q^3), \\
\cH_1(\tau, z)=& \sum_{\substack{n,\, r\in \ZZ \\ 2\nmid n, \; 2\nmid r \\ 4 n\geq r^2}}24 H\left(4 n-r^2\right) q^n \zeta^r 
			=8q\zeta^{\pm 1} + 8(\zeta^{\pm 3}+3\zeta^{\pm 1}) q^3 + O(q^5). 
\end{align}
By Table \ref{tab:dim-4}, $\cH_0$ and $\cH_1$ span the vector space $J_{2,1}(\Gamma_0(4))$. By Lemma \ref{lem:U} and Lemma \ref{lem:V}, we construct the following holomorphic Jacobi forms of weight $2$ and index $4$ on $\Gamma_0(4)$:
\begin{align}
\Big(\cH_0|_{2,1}V_2\Big)(2\tau, 2z) &= \sum_{\substack{n,\, r\in \ZZ \\ 2|n,\; 2|r \\ 4n\geq r^2}} 12H^{(2)}(4n-r^2) q^n\zeta^{2r} = 1+6(\zeta^{\pm 4} + 2)q^2 + O(q^4),\\
\cH_0(\tau,2z)&=\sum_{\substack{n,\, r\in \ZZ \\ 2|r \\ 4n\geq r^2}} 12H^{(2)}(4n-r^2) q^{n}\zeta^{2r} = 1+(\zeta^{\pm 4} + 6)q + O(q^2),\\
\cH_1(\tau,2z)&=\sum_{\substack{n,\, r\in \ZZ \\ 2\nmid n,\; 2\nmid r \\ 4n\geq r^2}} 24H(4n-r^2) q^{n}\zeta^{2r} = 8q\zeta^{\pm 2} + O(q^3).
\end{align}
Clearly, these forms are linearly independent and thus span $J_{2,4}(\Gamma_0(4))$. By comparing Fourier coefficients, we find
\begin{equation}
\vartheta(2\tau,2z)^4= \cH_0(\tau,2z) - \frac{1}{2}\cH_1(\tau,2z) -  \Big(\cH_0|_{2,1}V_2\Big)(2\tau, 2z). 
\end{equation}
This yields the desired formula in Theorem \ref{th:power-4}. 

\subsection{A proof of Theorem \ref{th:4-polygonal}}
Let $m$ be a positive integer. By \eqref{eq:theta}, we have
$$
\vartheta\big(m\tau,1/2-\tau\big) = i q^{m/8-1/2}\sum_{\ell \in \ZZ} q^{p_{m+2}(\ell)}, \quad p_{m+2}(\ell) = \frac{m\ell^2-(m-2)\ell}{2}. 
$$
This proves the generating function \eqref{eq:generating}. We then derive the desired formula for $R_{m+2}(N)$ in Theorem \ref{th:4-polygonal} by taking $\tau=m\tau'/2$ and $z=1/4-\tau'/2$ in \eqref{eq:power-4} of Theorem \ref{th:power-4}. 

\subsection{A proof of Corollary \ref{cor:4-polygonal}}
If the second sum in \eqref{eq:R4} is non-empty, then $4n-r^2$ never vanishes and thus $H(4n-r^2)>0$, since $r$ is odd. Recall that $H(D)\geq 0$ for any $D>0$ and $H^{(2)}(D)\geq 0$ for any $D\geq 0$. Therefore, every single term on the right-hand side of \eqref{eq:R4} is non-negative. Assume $m$ is odd. Then there exists an odd integer $r_0$ with $|r_0|\leq m$ such that $n_0=2(N-2+r)/m+1$ is also an odd integer. When $N\geq m^3/8+m/2+2$, $4n_0-r_0^2$ is always positive, which follows that $R_{m+2,4}(N)\geq 12H(4n_0-r_0^2)>0$. In a similar way, we can prove the corollary for even $m$ and odd $N$. 

\subsection{A proof of Corollary \ref{cor:r=1}}
When $m=6$ and $N$ is odd, $R_{8,4}(N)$ can be simplified as 
$$
R_{8,4}(N) = \sum_{\substack{3n-r-1=N \\ 2\nmid n, \,  2\nmid r  \\ 4 n-r^2 \geqslant 0}}12 H(4n-r^2).
$$
There exist different odd integers $r_0$, $r_1$, $r_2$ such that $|r_0|\leq|r_1|<|r_2|\leq9$ and $n_i=(N+r_i+1)/3$ are odd integers for $i=0,1,2$. Indeed, $(r_0,r_1,r_2)$ can be chosen as $(1,-5,7)$, $(-1,5,-7)$ and $(-3,3,9)$ if $N$ is of the  form $6x+1$, $6x+3$ and $6x-1$, respectively. If $N>53$, then $4n_i-r_i^2>3$ for $i=1,2,3$. Clearly, $4n_i-r_i^2$ is of the form $8y+3$. It follows that 
$$
R_{8,4}(N)\geq \sum_{i=1}^3 12H(4n_i-r_i^2) >24.
$$
Hence $r(N)>1$ for odd $N> 53$. We then prove the corollary by computer calculations.

\subsection{A proof of Proposition \ref{prop:eta-4}} 
It is easy to check that 
$$
\vartheta(3\tau,2\tau)=-q^{-5/8}\prod_{n=1}^\infty(1-q^n).
$$
We then establish the desired identity by taking $\tau=3\tau'/2$ and $z=\tau'$ in \eqref{eq:power-4}. 

\subsection{A proof of Proposition \ref{prop:eta-12}}
It is well-known that $\vartheta(\tau,0)=0$ and $\frac{\partial \vartheta(\tau,z)}{\partial z}\Big|_{z=0}=2\pi i \eta(\tau)^3$, where $\eta(\tau)=q^{1/24}\prod_{n=1}^\infty (1-q^n)$. It follows that 
$$
\dfrac{\partial^4 \vartheta(2\tau, 2z)^4}{\partial z^4}\Big|_{z=0}=6144\pi^4 \eta(2\tau)^{12}. 
$$
We then deduce Proposition \ref{prop:eta-12} from \eqref{eq:power-4} and the above identity.

\subsection{A proof of Proposition \ref{prop:eta-1-2-8}}
We consider the space of holomorphic Jacobi forms of weight $2$ and index $9$ on $\Gamma_0(2)$. By Table \ref{tab:dim-2}, this vector space is of dimension $2$. Moreover, it can be generated by $\cH(\tau,3z)$ and $\big(\cH|_{2,1}V_9\big)(\tau,z)$.  It is easy to check that
$$
\psi_{2,9}(\tau,z):=\frac{\vartheta(\tau,z)^2\vartheta(\tau,2z)\vartheta(2\tau,2z)^2\vartheta(2\tau,4z)}{\eta(\tau)\eta(2\tau)} \in J_{2,9}(\Gamma_0(2)).
$$
By comparing their Fourier expansions, we find  
\begin{equation}
\begin{split}
\psi_{2,9}(\tau,z)&= \frac{13}{12}\cH(\tau,3z)-\frac{1}{12}\Big(\cH|_{2,1}V_9\Big)(\tau,z)\\
&=13 \sum_{\substack{n,\, r\in \ZZ, \; 3|r \\ 36n\geq r^{2} }} H^{(2)}\left(\frac{36n-r^2}{9}\right)  q^n\zeta^r - \sum_{\substack{n,\, r\in \ZZ \\ 36n\geq r^{2}}} \sum_{\substack{a>0 \\ a \mid (n,r,9) }} a H^{(2)} \left(\frac{36n-r^2}{a^2}\right)q^n\zeta^r. 
\end{split}    
\end{equation}
Combining this formula and the identity
$$
\dfrac{\partial^6 \psi_{2,9}(\tau,z)}{\partial z^6}\Big|_{z=0}=-1474560\pi^6 \eta(\tau)^8\eta(2\tau)^{8},
$$
we prove Proposition \ref{prop:eta-1-2-8}.

\section{Proofs of Theorem \ref{th:four-squares} and its corollary} \label{sec:4-square}
For any $u=(x_1,x_2,x_3,x_4)\in \ZZ^4$ and $v=(y_1,y_2,y_3,y_4)\in\ZZ^4$, we define 
$$
\latt{u,v}=x_1y_1+x_2y_2+x_3y_3+x_4y_4. 
$$
The representative number introduced in Theorem \ref{th:four-squares} is also given by
\begin{equation}
\rho(n, r, m)=\#\left\{ u,\, v \in \ZZ^4 : \latt{u,u}=n, \; \latt{v,v}=m,\; \latt{v,u}=r \right\}.   
\end{equation}
We then derive the generating function
\begin{equation}
\begin{split}
\Theta_{4A_1}(Z):=& \sum_{\substack{n,\, m\in \NN,\; r\in \ZZ\\ nm\geq r^2}} \rho(n,r,m) \exp\Big( 2\pi i\big(n\tau+2rz+m\tau'\big) \Big) \\
=& \sum_{u,\,v\in \ZZ^4} \exp\Big( 2\pi i\big(\latt{u,u}\tau+2\latt{u,v}z+\latt{v,v}\tau'\big) \Big),
\end{split}   
\end{equation}
where $Z=\begin{psmallmatrix}
\tau & z \\ z & \tau'    
\end{psmallmatrix}$ lies in $\HH_2$, namely the Siegel upper half space of degree two. Since the even positive-definite lattice $4A_1$, that is, $\ZZ^4$ with bilinear form $2\latt{-,-}$, is of level $4$ and discriminant $2^4$, we conclude from \cite{And87} that $\Theta_{4A_1}(Z)$ is a Siegel modular form of weight $2$ and trivial character on the level $4$ congruence subgroup
$$
\Gamma_0^2(4)=\big\{ M \in \Sp_4(\ZZ) \; : \; M \equiv \begin{psmallmatrix}
* & * \\ 0 & *     
\end{psmallmatrix} \; (\mathrm{mod}\; 4) \big\}.
$$
Let $M_2(\Gamma_0^2(4))$ denote the space of Siegel modular forms of weight $2$ and trivial character on $\Gamma_0^2(4)$. It is known that $M_2(\Gamma_0^2(4))$ is of dimension $3$ (see e.g. \cite[Appendix B.2]{RSY23}). We now construct a basis of this space whose Fourier expansions involve the Hurwitz class numbers. 

By \cite[Theorem 2.2]{CG11}, we have the following Siegel modular forms of weight $2$ on $\Gamma_0^2(4)$:
\begin{align*}
F_0(Z) &= \sum_{m=0}^{\infty} \big(\cH_0 |V_m\big)(\tau,z)e^{2\pi i m\tau'}=\frac{1}{24}+\sum_{n=1}^\infty\sum_{\substack{a|n \\ 2\nmid a}}aq^n + \cH_0(\tau,z)e^{2\pi i\tau'}+O(e^{4\pi i\tau'}),\\
F_1(Z) &= \sum_{m=1}^{\infty} \big(\cH_1 |V_m\big)(\tau,z)e^{2\pi i m\tau'}=\cH_1(\tau,z)e^{2\pi i\tau'}+O(e^{4\pi i\tau'}),\\
G(Z) &= \sum_{m=0}^{\infty} \big(\cH |V_m\big)(2\tau,2z)e^{4\pi i m\tau'}=\frac{1}{24}+\sum_{n=1}^\infty\sum_{\substack{a|n \\ 2\nmid a}}aq^{2n} + \cH(2\tau,2z)e^{4\pi i\tau'}+O(e^{8\pi i\tau'}).
\end{align*}
By comparing Fourier coefficients, $F_0$, $F_1$ and $G$ are linearly independent and thus form a basis of $M_2(\Gamma_0^2(4))$. Moreover, their Fourier expansions can be expressed as linear combinations of $H^{(2)}(-)$ and $H(-)$ by Lemma \ref{lem:V}.  This is sufficient to establish the equality
\begin{equation}
\Theta_{4A_1}(Z)=8F_0(Z)+16G(Z)    
\end{equation}
and further deduce the desired formula in Theorem \ref{th:four-squares}. 

We now prove Lemma \ref{lem:positive}.  Let us write $D=2^{4k+l}m$ with $k,m\in\NN$, $m$ odd and $l\in \{0,1,2,3\}$. Gauss--Legendre's three-square theorem asserts that every non-negative integer $N$ is the sum of three squares if and only if $N$ is not of the form $4^a(8b+7)$ with $a,b\in\NN$.  By \eqref{eq:Gauss}, if $s\in\NN$ satisfies that $D-16^s$ is non-negative and not of the form $4^a(8b+7)$, then $H^{(2)}\big(4(D-16^s)\big)>0$. Let us set $s=k+c$ and
$$
c=\left\{\begin{array}{ll}1, & l=0,\, m\equiv 1\,(\bmod \, 4) \\ 0, &  l=0,\, m \equiv 3 \,(\bmod \, 4) \\ 0, & l=1,2 \\ 1, & l=3.
\end{array}\right.
$$
Then $D-16^s=4^{2k}(2^lm-16^c)$ and there is an integer $b$ such that
$$
2^{l}m-16^c =\left\{\begin{array}{ll}4b+1, & l=0,\, m\equiv 1\,(\bmod \, 4) \\ 4b+2, &  l=0,\, m \equiv 3 \,(\bmod \, 4) \\ 4b+1, & l=1 \\ 8b+3, & l=2 \\ 4(4b+2), & l=3. 
\end{array}\right.
$$
In the above, $2^{l}m-16^c\geq0$ holds except for 
\begin{enumerate}
\item $l=0$ and $m=1,5,9,13$;
\item $l=3$ and $m=1$.
\end{enumerate}
In the first case, we choose $c=0$ and thus $D-16^s=4^{2k}(m-1)\geq 0$ is not of the form $4^a(8b+7)$.  In the second case, if $D-16^s=16^k\cdot 8 - 16^s$ is positive, then $k\geq s$ and therefore $D-16^s$ is of the form $4^a(8b+7)$, which follows that $H^{(2)}\big(4(D-16^s)\big)=0$ for any $s\in\NN$. This yields the first part of Lemma \ref{lem:positive}. From
$$
\rho(n,4^s,m)\geq 96H^{(2)}\big(4(nm-16^s)\big)
$$
we deduce the second part of the lemma. 

We finally prove Corollary \ref{coro:sun}. Recall from the introduction that $\Orth(\ZZ^4)$ acts on $\mathcal{S}_m$ and $t_m$ denotes the number of orbits.  Notice that $R_{4,4}(4N)=R_{4,4}(N)$ for any even positive integer $N$. Clearly, if $R_{4,4}(m)>|\Orth(\ZZ^4)|=384$ then $t_m>1$. By combining these facts and computer calculations, we find that $t_m=1$ if and only if $m=1$, $3$, $5$, $7$, $11$, $15$, $23$, or $m=2\cdot 4^t$, $6\cdot 4^t$, $14\cdot 4^t$ for $t\in\NN$. When $m=9$, $t_m=2$ and the two orbits are represented by $v_1=(2,2,1,0)$ and $v_2=(3,0,0,0)$. However, $\latt{v_2,u}$ is never of the form $4^s$ for any $u\in \ZZ^4$. It follows that $\rho(n,4^s,9)=\delta_1\cdot \rho_{v_1}(n,4^s)$. We therefore prove the corollary by these facts and Lemma \ref{lem:positive}. 

\section{Appendix}
In this appendix we use Jacobi forms on $\Gamma_0(3)$ to derive formulas involving $H(D)$ for the Fourier coefficients of $\vartheta(\tau,z)^2\vartheta(3\tau,3z)^2$. As applications, we obtain formulas for the representative number of positive integers as sums of four generalized $m$-gonal numbers with coefficients $(1,1,3,3)$, and for the Fourier coefficients of cusp forms $\eta(\tau)^2\eta(3\tau)^2$ and $\eta(\tau)^6\eta(3\tau)^6$. The proof is pretty much the same as before, so we just sketch it. 

We verify that 
$$
\vartheta(\tau,z)^2\vartheta(3\tau,3z)^2=(\zeta^{\pm4}-2\zeta^{\pm3}+\zeta^{\pm2}-2\zeta^{\pm 1}+4)q+O(q^2)
$$
is a holomorphic Jacobi form of weight $2$ and index $4$ with trivial character on $\Gamma_0(3)$. By \cite{BSZ19},  
\begin{equation}
\cH^\ast(\tau,z)= \sum_{\substack{ n\in\NN,\, r\in\ZZ \\ 4n-r^{2}\geq 0}} H^{(3)}\left(4 n-r^2\right) q^n \zeta^r
\end{equation}
is a basis of $J_{2,1}(\Gamma_0(3))$, where  
$$
H^{(3)}(D)=H(9D)-3H(D)
$$
if $D\equiv 0,3\; (\mathrm{mod}\; 4)$, and $H^{(3)}(D)=0$ otherwise. We find that $\dim J_{2,4}(\Gamma_0(3))=2$ by the algorithm of Section \ref{sec:Jacobi}.  Lemmas \ref{lem:U}--\ref{lem:V} yield that $\cH^\ast(\tau,2z)$ and $\big(\cH^\ast|_{2,1}V_4\big)$ span $J_{2,4}(\Gamma_0(3))$. We then derive
\begin{equation}\label{eq:1133}
\begin{split}
\vartheta(\tau,z)^2\vartheta(3\tau,3z)^2=&7\cH^\ast(\tau,2z)-\big(\cH^\ast|V_4\big)(\tau,z) \\
=&\sum_{\substack{n\in\NN,\, r\in\ZZ \\ 4n-r^{2} \geq 0}} 7H^{(3)}(4n-r^2) q^n\zeta^{2r} - \sum_{\substack{n\in\NN,\, r\in\ZZ \\ 16n-r^{2} \geq 0}} H^{(3)}(16n-r^2)q^n\zeta^r\\
&-\sum_{\substack{n\in\NN,\, r\in\ZZ \\ 8n-r^{2} \geq 0}} 2H^{(3)}(8n-r^2)q^{2n}\zeta^{2r}-\sum_{\substack{n\in\NN,\, r\in\ZZ \\ 4n-r^{2} \geq 0}}4H^{(3)}(4n-r^2)q^{4n}\zeta^{4r}.
\end{split}    
\end{equation}

For any $N\in\NN$ we define the representative number
\begin{equation}
R_{m}^{\ast}(N)=\# \{ (\ell_1,\ell_2,\ell_3,\ell_4) \in \ZZ^4 : N = p_{m}(\ell_1)+ p_{m}(\ell_2)+3p_{m}(\ell_3)+3 p_{m}(\ell_4)\}.     
\end{equation}
From the generating function 
$$\vartheta\big(m\tau,1/2-\tau\big)^2\vartheta\big(3m\tau,3/2-3\tau\big)^2=q^{m-4}\sum_{N=0}^{\infty}R_{m+2}^{\ast}(N)q^N
$$
and \eqref{eq:1133} we deduce
\begin{equation}
\begin{split}
R_{m+2}^{\ast}(N)=&\sum_{\substack{n\in\NN,\, r\in\ZZ \\ mn-2r+4-m=N\\4n-r^2\geq0}}7H^{(3)}(4n-r^2)-\sum_{\substack{n\in\NN,\, r\in\ZZ \\ mn-r+4-m=N\\16n-r^2\geq0}}(-1)^rH^{(3)}(16n-r^2) \\
&-\sum_{\substack{n\in\NN,\, r\in\ZZ \\ 2mn-2r+4-m=N\\8n-r^2\geq0}}2H^{(3)}(8n-r^2)-\sum_{\substack{n\in\NN,\, r\in\ZZ \\ 4mn-4r+4-m=N\\4n-r^2\geq0}}4H^{(3)}(4n-r^2).   
\end{split}    
\end{equation}

Let us write $\prod_{n=1}^\infty(1-q^n)^2(1-q^{3n})^2=\sum_{n=0}^\infty a_nq^n$. From the identity
$$
\vartheta(3\tau,2\tau)^2\vartheta(9\tau,6\tau)^2=q^{-16/3}\eta(\tau)^2\eta(3\tau)^2
$$
and \eqref{eq:1133} we derive
\begin{equation}
\begin{split}
a_N=&\sum_{\substack{n\in\NN,\, r\in\ZZ \\ 3n+4r+5=N \\4n-r^2\geq0}} 7H^{(3)}(4n-r^2)-\sum_{\substack{n\in\NN,\, r\in\ZZ \\ 3n+2r+5=N \\16n-r^2\geq0}} H^{(3)}(16n-r^2)\\
&-\sum_{\substack{n\in\NN,\, r\in\ZZ \\ 6n+4r+5=N \\8n-r^2\geq0}}2H^{(3)}(8n-r^2)-\sum_{\substack{n\in\NN,\, r\in\ZZ \\ 12n+8r+5=N \\4n-r^2\geq0}}4H^{(3)}(4n-r^2). 
\end{split}
\end{equation}

From the identity
$$
\frac{\partial^4 \vartheta(\tau,z)^2\vartheta(3\tau,3z)^2}{\partial z^4} \Big|_{z=0}=3456\pi^4\eta(\tau)^6\eta(3\tau)^6
$$
and \eqref{eq:1133} we conclude that $\eta(\tau)^6\eta(3\tau)^6$ has the Fourier expansion
\begin{equation}
\begin{split}
&\frac{14}{27}\sum_{n=1}^\infty  \sum_{\substack{r\in\ZZ \\ r^{2}\leq 4n}}r^4H^{(3)}(4n-r^2)q^{n} - \frac{1}{216}\sum_{n=1}^\infty \sum_{\substack{r\in\ZZ \\ r^{2}\leq 16n}}r^4H^{(3)}(16n-r^2)q^{n}\\
-&\frac{4}{27}\sum_{n=1}^\infty\sum_{\substack{r\in\ZZ \\ r^{2}\leq 8n}}r^4H^{(3)}(8n-r^2)q^{2n}-\frac{32}{27}\sum_{n=1}^\infty\sum_{\substack{r\in\ZZ \\ r^{2}\leq 4n}}r^4H^{(3)}(4n-r^2)q^{4n}. 
\end{split}
\end{equation}

\bibliographystyle{plainnat}
\bibliofont
\bibliography{refs}

\end{document}